
\magnification=\magstep1

\input amstex

\documentstyle{amsppt}

\TagsOnRight

\def\cite#1{{\rm [#1]}}
\def\ol{\overline }

\leftheadtext{B. AUPETIT, E. MAKAI, Jr., M. MBEKHTA, and J. ZEM\'ANEK}

\rightheadtext{LIFTINGS OF ANALYTIC FAMILIES OF IDEMPOTENTS}

\topmatter

\title
{\centerline{Local and global liftings of analytic families}}
{\centerline{of idempotents in Banach algebras}}
\vskip.2cm
{\centerline{{\rm{Dedicated to the 100th anniversary}}}} 
{\centerline{{\rm{of Professor B\'ela Sz\H okefalvi-Nagy}}}}
\endtitle

\author
{Bernard Aupetit, Endre Makai, Jr., Mostafa Mbekhta, and Jaroslav Zem\'anek}
\endauthor

\address
Bernard Aupetit, 
\newline
\indent D\'epartement de Math\'ematiques, 
\newline
\indent Universit\'e Laval, 
\newline
\indent Qu\'ebec, G1K 7P4, 
\newline
\indent CANADA
\endaddress

\address
Endre Makai, Jr.,
\newline
\indent Alfr\'ed R\'enyi Mathematical Institute,
\newline
\indent Hungarian Academy of Sciences,
\newline
\indent H-1364 Budapest, P.O. Box 127,
\newline
\indent HUNGARY
\endaddress
\email
makai.endre\@renyi.mta.hu, {\rm{http://www.renyi.mta.hu/\~{}makai}}
\endemail

\address
Mostafa Mbekhta,
\newline
\indent Laboratoire P. Painlev\'e, UMR-CNRS 8524
\newline
\indent U.F.R. de Math\'ematiques Pures et Appliqu\'ees, 
B\^at. M2,
\newline
\indent Universit\'e des Sciences et Technologies de Lille,
\newline
\indent F-59655 Villeneuve d'Ascq,
\newline
\indent FRANCE
\endaddress
\email
mostafa.mbekhta\@math.univ-lille1.fr
\endemail

\address
Jaroslav Zem\'anek,
\newline
\indent Institute of Mathematics, Polish Academy of Sciences,
\newline
\indent 00-956 Warsaw, P.O. Box 21,
\newline
\indent POLAND
\endaddress
\email
zemanek\@impan.pl
\endemail

\thanks Research (partially) supported by Hungarian National Foundation for
Scientific Research, grant numbers K68398, K75016, K81146, 
and the European Commission Project TODEQ (MTKD-CT-2005-030042).
\endthanks



\abstract
Generalizing results of our earlier paper, we investigate the following 
question.
Let $\pi(\lambda) : A \to B$ be an analytic family of surjective 
homomorphisms between two Banach algebras, and $q(\lambda)$ an analytic 
family of idempotents in~$B$.
We want to find an analytic family $p(\lambda)$ of idempotents in~$A$, 
lifting $q(\lambda)$, i.e., such that $\pi(\lambda)p(\lambda) = q(\lambda)$, 
under hypotheses of the type that the elements of\/ $\text{\rm Ker}\, 
\pi(\lambda)$ have small spectra.
For spectra which do not disconnect $\Bbb C$ we obtain a local lifting theorem.
For real analytic families of surjective $^*$-homomorphisms (for continuous
involutions) and self-adjoint 
idempotents we obtain a local lifting theorem, for totally disconnected 
spectra.
We obtain a global lifting theorem if the spectra of the elements in 
$\text{\rm Ker}\, \pi(\lambda)$ are $\{0\}$, both in the analytic case, 
and, for $^*$-algebras (with continuous involutions)
and self-adjoint idempotents, in the real analytic case.
Here even an at most countably infinite set of mutually orthogonal analytic 
families of idempotents can be lifted to mutually orthogonal analytic 
families of idempotents. In the proofs, spectral theory is combined with
complex analysis and general topology, and even a connection with potential
theory is mentioned.

\endabstract

\endtopmatter

\document


\heading
1. Notations
\endheading

For a Banach algebra $A$, let
$$
E(A) := \{ p \in A: \ p^2 = p\}
$$
be the set of idempotents.
If there is an involution $^*$ on $A$, let
$$
S(A) := \{ p \in A: \ p^2 = p = p^*\}
$$
be the set of self-adjoint idempotents.

For a Banach space $X$, let $\Cal B(X)$ or $\Cal K(X)$ denote the Banach 
algebras of all bounded or all compact linear operators on~$X$, respectively.
(For two Banach spaces $X$, $Y$, let $\Cal B(X, Y)$ denote the Banach space 
of all bounded linear operators $X \to Y$.)
The quotient $\Cal C(X) := \Cal B(X) / \Cal K(X)$ is called the 
{\it Calkin algebra} (as it was first studied in \cite{C} for the 
Hilbert space case).
Let
$$
\pi : \ \Cal B(X) \to \Cal C(X)
$$
be the canonical map.
More generally, we will consider situations where $\pi : A \to B$ is a 
surjective map between Banach algebras.

In this paper, we will study the relationship between $E(A)$ and $E(A/I)$, 
where $I$ is a closed two-sided ideal in~$A$.
Observe that
$$
\pi E (A) \subset E(A/I),
$$
where $\pi : A \to A/I$ is the canonical map.
Thus, the lifting problem consists in establishing the converse inclusion.
See Gramsch \cite{Gr} for a case, when the converse inclusion does not hold 
(cf.\ also \cite{AMMZ, p.~26}).

In this paper, the linear structures are considered over the complex scalars 
$\Bbb C$, the algebras are assumed to be unital (unless explicitly stated 
contrarily), and all Banach algebra 
homomorphisms are assumed to be continuous, and preserving 1.
We denote by $\sigma(x)$ the spectrum of an element $x$ of a Banach algebra.
When considering spectra of elements of various kernels of mappings $A \to B$,
which kernels are therefore subsets of the Banach algebra $A$, the spectra are
considered with respect to $A$.


\heading
2. Lifting of idempotents and local lifting of analytic families of idempotents
\endheading

The lifting of a single self-adjoint idempotent was first studied for 
the case $A = \Cal B(H)$ and $I = \Cal K(H)$, where $H$ is a Hilbert space, 
with the result that $\pi S(\Cal B(H)) = S(\Cal C(H))$, see 
\cite{C, Theorem 2.4} and \cite{dlH, Proposition 7}.
For a general Banach algebra $A$ and $I = \text{\rm rad}\, A$, the 
lifting property $\pi E(A) = E(A / \text{\rm rad}\, A)$ can be found in 
\cite{Ri, p. 58, Theorem 2.3.9} and \cite{Ka, p.~125}.
The formula $\pi E(\Cal B(H)) = E(\Cal C(H))$ was obtained in 
\cite{La, Theorem 2}.
Labrousse actually proved that each analytic family of 
idempotents of\/ $E(\Cal C(H))$ can be locally lifted to an analytic family 
of idempotents in $E(\Cal B(H))$.
That is, for each analytic map $q : U \to E(\Cal C(H))$, where $U$ is an 
open set in $\Bbb C$, any point of\/ $U$ has an open neighbourhood 
$V \subset U$, such that there exists an analytic map 
$p : V \to E(\Cal B(H))$ with $\pi p(\lambda) = q(\lambda)$ for each 
$\lambda \in V$, with $\pi : \Cal B(H) \to \Cal C(H)$ being the canonical map.
Recall that every compact operator on a Banach space has a spectrum that is 
finite, or is a sequence converging to~$0$, together with $0$ 
(\cite{DS, p. 579, Theorem VII. 4.5}).
Such a compact set is {\it{totally disconnected}} (i.e., it contains
no connected subset consisting of 
more than one points), which implies, due to 
compactness, that its complement in $\Bbb C$ is connected 
(\cite{Ku, p. 466, Section 59, II, Theorem~1, and p. 188,
Section 47, VIII, Theorem~1}). (Cf. also the remarks before the proof of
our Theorem 2.)
This raises the question, what can be said in general about the equality 
$\pi E(A) = E(A/I)$, if we make restrictions on the spectra of the elements 
in $I = \pi^{-1}(0) = \text{\rm Ker}\, \pi$.

In \cite{AMMZ} there were proved four theorems of this type.

\newpage

\proclaim{Theorem A} {\rm (\cite{AMMZ, Theorem 3.1 and Corollary 3.3})}
Let $\pi : A \to B$ be a surjective homomorphism between Banach algebras.
Suppose that the spectrum of each element of\/ $\text{\rm Ker}\, \pi$ is 
totally disconnected.
Let $U \subset \Bbb C$ be open, and let $q : U \to E(B)$ be an analytic map.
Then any point of\/ $U$ has an open neighbourhood $V \subset U$, such that 
there exists an analytic map $p : V \to E(A)$ with 
$\pi p(\lambda) = q(\lambda)$ for each $\lambda \in V$.
In particular, for $q$ constant, we obtain $\pi E(A) = E(B)$.
\endproclaim

\proclaim{Theorem B} {\rm (\cite{AMMZ, Corollary 3.3})}
Assume the hypotheses of Theorem~A, and let $A$ and $B$ possess continuous 
involutions, and let $\pi$ be a $^*$-homomorphism.
Then $\pi S(A) = S(B)$.
\endproclaim

\proclaim{Theorem C} {\rm (\cite{AMMZ, p.~26, paragraph following the proof of
Corollary 3.3})}
Let $\pi : A \to B$ be a surjective homomorphism of Banach algebras.
Suppose that the spectrum of each element of\/ $\text{\rm Ker}\, \pi$ does 
not disconnect~$\Bbb C$.
Then $\pi E (A) = E(B)$.
\endproclaim

\proclaim{Theorem D} {\rm (\cite{AMMZ, Corollary 3.2})}
For the canonical map $\pi : A \to A / \text{\rm rad}\, A$ we may choose 
$V = U$ in Theorem~A, i.e., there is a global lifting of~$q(\lambda)$.
\endproclaim

We cite below a restriction of a theorem from \cite{Ka, pp.\ 125--126} to 
Banach algebras.
Before this we recall that two idempotents $e$ and $f$ of a Banach algebra 
are {\it orthogonal}, if\/ $ef = fe = 0$.

\proclaim{Theorem E} {\rm (\cite{Ka, pp.\ 125--126})}
If\/ $A$ is a Banach algebra and $q_1, q_2, \dots$ is a finite or countably 
infinite set of orthogonal idempotents in $A / \text{\rm rad}\, A$, then 
there are orthogonal idempotents $p_1, p_2, \dots$ in~$A$ (indexed by the 
same index set as the $q_i$'s), such that $\pi(p_i) = q_i$, for all $i$'s,
where $\pi$ is 
the canonical map $A \to A / \text{\rm rad}\, A$.
\endproclaim

In this paper we generalize the above theorems.

For a more comprehensive literature on the structure of the set of 
idempotents, we refer to \cite{AMMZ, \S2, Historical background}.


\heading
3. Local and global liftings of analytic families of idempotents, for analytic 
families of surjective homomorphisms
\endheading

We give a common generalization of Theorems~A and C.

\proclaim{Theorem 1}
Let $U$ be an open subset of\/ $\Bbb C$, containing $0$.
Let $A$ and $B$ be Banach algebras, and let $\pi : U \to \Cal B(A, B)$ be an 
analytic map, whose values are homomorphisms $A \to B$, such that $\pi(0)$ 
is surjective.
Suppose that the spectrum of each element of\/ $\text{\rm Ker}\, \pi(0)$ 
does not disconnect $\Bbb C$.
Let $q: U \to E(B)$ be an analytic map.
Then there exist an open set $V \subset \Bbb C$, such that 
$0 \in V \subset U$, and an analytic map $p : V \to E(A)$, such that 
$\pi(\lambda) p(\lambda) = q(\lambda)$ for each $\lambda \in V$.
\endproclaim

We prove a generalization of Theorem~B.
By a {\it real analytic map} from an open subset $G$ of\/ $\Bbb R$ to a 
Banach space, we mean a map $f$ that for each $x_0 \in G$ has locally a 
power series expansion $f(x) = \sum\limits^\infty_0 a_n(x - x_0)^n$.
In the real analytic case, we need a stronger spectral assumption (total
disconnectedness), 
again, for one of the kernels (namely, at $0$), merely. Then we have a variant
of Theorem 1, for Banach algebras with continuous involutions.

\proclaim{Theorem 2}
Let $G$ be an open subset of~$\Bbb R$, containing~$0$.
Let $A$ and $B$ be Banach 
\newpage

algebras with continuous involutions, and let 
$\pi : G \to \Cal B(A, B)$ be a real analytic map, whose values are 
$^*$-homomorphisms $A \to B$, such that $\pi(0)$ is surjective.
Suppose that the spectrum of each element of\/ $\text{\rm Ker}\, \pi(0)$ is 
totally disconnected.
Let $q : G \to S(B)$ be a real analytic map.
Then there exist an open set $H \subset \Bbb R$, such that 
$0 \in H \subset G$, and a real analytic map $p : H \to S(A)$, such that 
$\pi(\lambda) p(\lambda) = q(\lambda)$ for each $\lambda \in H$.
\endproclaim

We give a common generalization of Theorems~D and~E.
Recall that for the canonical map $\pi : A \to A / \text{\rm rad}\, A$ we 
have $\sigma(a) = \sigma(\pi(a))$ for each $a \in A$, cf. [A79, p. 2, Lemme
1.1.2], or
\cite{A91, p. 35, Theorem
3.1.5}. Here the strongest spectral assumption is imposed on all the kernels
in question. 

\proclaim{Theorem 3}
Let $U$ be an open subset of\/ $\Bbb C$.
Let $A$ and $B$ be Banach algebras, and let $\pi : U \to \Cal B(A, B)$ be an 
analytic map, whose values are surjective homomorphisms $A \to B$.
Suppose that the spectrum of each element 
of\/ $\text{\rm Ker}\, \pi(\lambda)$, for each $\lambda \in U$, is $\{0\}$.
Let $q_1, q_2, \dots : U\to E(B)$ be finitely or countably infinitely many 
analytic maps, where for each $i \neq j$ and each $\lambda \in U$ we have 
that $q_i(\lambda)$ and $q_j(\lambda)$ are orthogonal.
Then there exist analytic maps $p_1, p_2, \dots : U \to E(A)$ 
(indexed by the same index set as the $q_i$'s), where for each $i \neq j$ 
and each $\lambda \in U$ we have that $p_i(\lambda)$ and $p_j(\lambda)$ 
are orthogonal, and such that $\pi(\lambda) p_i(\lambda) = q_i(\lambda)$ 
holds for each index~$i$ and each $\lambda \in U$.
\endproclaim

We give a variant of Theorem~3 for Banach algebras with continuous involutions.

\proclaim{Theorem 4}
Let $G$ be an open subset of~$\Bbb R$.
Let $A$ and $B$ be Banach algebras with continuous involutions, and let 
$\pi : G \to \Cal B(A, B)$ be a real analytic map, whose values are 
surjective $^*$-homomorphisms $A \to B$.
Suppose that the spectrum of each element of\/ $\text{\rm Ker}\, 
\pi(\lambda)$, for each $\lambda \in G$, is~$\{0\}$.
Let $q_1, q_2, \dots : G \to S(B)$ be finitely or countably infinitely many 
real analytic maps, where for each $i \neq j$ and each $\lambda \in G$ we 
have that $q_i(\lambda)$ and $q_j(\lambda)$ are orthogonal.
Then there exist real analytic maps $p_1, p_2, \dots : G \to S(A)$ 
(indexed by the same index set as the $q_i$'s), where for each $i \neq j$ 
and each $\lambda \in G$ we have that $p_i(\lambda)$ and $p_j(\lambda)$ 
are orthogonal, and such that $\pi(\lambda) p_i(\lambda) = q_i(\lambda)$ 
holds for each index~$i$ and each $\lambda \in G$.
\endproclaim

In the next two theorems, the strongest spectral assumption is imposed only on
one of the kernels (namely, at $0$). But then we are able to lift only a
finite set of analytic families of orthogonal idempotents, and only locally.

\proclaim{Theorem 5}
Let $U$ be an open subset of\/ $\Bbb C$, containing $0$.
Let $A$ and $B$ be Banach algebras, and let $\pi : U \to \Cal B(A, B)$ be an 
analytic map, whose values are homomorphisms $A \to B$, such that $\pi (0)$ is
surjective. Suppose that the spectrum of each element 
of\/ $\text{\rm Ker}\, \pi(0)$ is $\{0\}$.
Let $q_1, \dots , q_n : U\to E(B)$ be finitely many 
analytic maps, where for each $i \neq j$ and each $\lambda \in U$ we have 
that $q_i(\lambda)$ and $q_j(\lambda)$ are orthogonal.
Then there exist an open set $V \subset {\Bbb C}$, such that $0 \in V \subset
U$, and 
analytic maps $p_1, \dots ,p_n: V \to E(A)$,
where for each $i \neq j$ 
and each $\lambda \in V$ we have that $p_i(\lambda)$ and $p_j(\lambda)$ 
are orthogonal, and such that $\pi(\lambda) p_i(\lambda) = q_i(\lambda)$ 
holds for each index~$i \in \{ 1, \dots , n \} $ and each $\lambda \in V$.
\endproclaim

The following Theorem 6 is a variant of Theorem 5, for Banach algebras with
continuous involutions.

\newpage

\proclaim{Theorem 6}
Let $G$ be an open subset of~$\Bbb R$, containing $0$. 
Let $A$ and $B$ be Banach algebras with continuous involutions, and let 
$\pi : G \to {\Cal B}(A, B)$ be a real analytic map, whose values are 
$^*$-homomorphisms $A \to B$, such that $\pi (0)$ is surjective. 
Suppose that the spectrum of each element of\/ $\text{\rm Ker}\, 
\pi(0)$ is~$\{0\}$.
Let $q_1, \dots , q_n: G \to S(B)$ be finitely many 
real analytic maps, where for each $i \neq j$ and each $\lambda \in G$ we 
have that $q_i(\lambda)$ and $q_j(\lambda)$ are orthogonal.
Then there exist an open set $H \subset {\Bbb R}$, such that $0 \in H \subset
G$, and
real analytic maps $p_1, \dots , p_n: H \to S(A)$, 
where for each $i \neq j$ 
and each $\lambda \in H$ we have that $p_i(\lambda)$ and $p_j(\lambda)$ 
are orthogonal, and such that $\pi(\lambda) p_i(\lambda) = q_i(\lambda)$ 
holds for each index~$i \in \{ 1, \dots , n \} $ and each $\lambda \in H$.
\endproclaim

{\bf{Remark 1.}}
In fact our theorems hold in greater generality, as the proofs given show this.
Theorems~1 and 5 hold for $0 \in U \subset \Bbb C^n$ open.
Theorems~2 and 6 hold for $0 \in G \subset \Bbb R^n$ open.
Theorem~3 holds if we replace $U$ by a Stein manifold.
Theorem~4 holds for $G \subset \Bbb R^n$ open, provided each connected 
component of\/ $U$ has a neighbourhood base consisting of domains of 
holomorphy, when we consider $\Bbb R^n$ as embedded in~$\Bbb C^n$ in the
canonical way.

\vskip.1cm

{\bf{Remark 2.}}
Another way of strengthening our theorems, for Theorems 3 and 4, is the
following. Retaining the hypothesis $U \subset {\Bbb C}$ or $G \subset {\Bbb
R}$ open, 
the spectral hypothesis, about the spectra of elements of Ker\,$\pi
(\lambda )$, is not necessary to be postulated 
for all $\lambda $ in $U$, or in $G$, respectively.
For Theorems 3 and 4 
it is sufficient to postulate it for some subsets of each connected component
of $U$, or $G$, respectively,
which have positive outer capacity (like, e.g., non-degenerate straight line
segments).
See the explanation in Remark 4, after the proof of Theorem 6. (For
capacity, and potential theory, cf. the books [HK], [A79], [A91], 
[Ra] and [AG].)

\vskip.1cm

Now we give a non-trivial example of (real) analytic families of surjective 
homomorphisms ($^*$-homomorphisms) between Banach algebras 
(Banach $^*$-algebras, with continuous involutions). Essentially the same
example was provided independently by Globevnik [Glo3] and 
Leiterer [Le2].

\vskip.1cm

{\bf{Example 1.}}
Let $D : = \{z \in \Bbb C : |z| < 1\}$, and let $n \geq 1$ be an integer.
Let $B$ be a Banach algebra (Banach $^*$-algebra, with a continuous 
involution); here the existence of unit is not required.
Let $A$ be the following Banach algebra (in general, also without unit)
of analytic functions $D^n \to B$:
$$
A \! :=\! \biggl\{\! f\! : D^n\! \to\! B\! :  
f\! =\! \sum^\infty_0 \! a_{k_1 \dots k_n} \, z_1^{k_1} \dots z_n^{k_n}, \, 
a_{k_1 \dots k_n} \! \in\! B, \, \|f\| := \! \sum^\infty_0\! 
\|a_{k_1 \dots k_n} \|\! < \!\infty \! \biggr\}\! .
$$
That is, $A$ is the set of absolutely convergent power series, when 
considered as functions $\ol D^n \to B$, with the norm given in the last
display formula.

If\/ $B$ is a Banach $^*$-algebra with continuous involution, we let
$$
f^* := \sum^\infty_0 a_{k_1 \dots k_n} ^* z_1^{k_1} \dots z_n^{k_n}.
$$
We observe that the now defined
involution on $A$ is continuous, as well. Indeed, let 
$$
\forall b \in B \,\,\,\,\,\, \| b^* \| \le C \cdot \| b \| \,,
$$ 

\newpage

for some $C \in (0, \infty )$, and let $C$ be
the smallest number such that this inequality holds.
Then, we have 
$$
\| f^* \| = \sum _0 ^{\infty } \| a_{k_1 \dots k_n } ^* \| \le C \cdot
\sum _0 ^{\infty } \| a_{k_1 \dots k_n } \| = C \cdot \| f \| \,,
$$
and also for this last inequality the number $C$ is the smallest number 
such that this inequality holds (this being true already for the smaller set
of all constant functions $f$ of norm $1$). So, the involution on $A$ is also
continuous (and the constant $C$ will not appear anymore in the sequel).

We define
$$
\pi(\lambda_1, \dots, \lambda_n) f:= f(\lambda_1, \dots, \lambda_n),
$$
for each $\lambda_1, \dots, \lambda_n \in D$.
This is a surjective homomorphism, and, in case of Banach $^*$-algebras, 
for $\lambda_1, \dots, \lambda_n \in (-1,1)$, this is even a $^*$-homomorphism.
We have the expansion
$$
\pi(\lambda_1, \dots, \lambda_n) f = \sum^\infty_0 a_{k_1 \dots k_n} 
\lambda_1^{k_1} \dots \lambda_n^{k_n} =: \sum^\infty_0 a_{k_1 \dots k_n} (f) 
\lambda_1^{k_1} \dots \lambda_n^{k_n}.
$$
Here the maps $f \mapsto a_{k_1 \dots k_n}(f)$ are linear operators of 
norm~$1$, so the above series gives the power series expansion 
of\/ $\pi(\lambda_1 , \dots, \lambda_n)$, for 
$\lambda_1, \dots, \lambda_n \in D$, hence 
$\pi(\lambda_1, \dots, \lambda_n)$ is, in fact, an analytic family on $D^n$.
It is interesting to observe, that $\pi(\lambda_1, \dots, \lambda_n)$,
for each $\lambda _1,..., \lambda _n \in D$,  
has norm $1$ --- hence their norm is constant --- which norm is
attained exactly for all constant functions $f(\lambda _1,...\lambda _n)=
a_{0...0}(f)$. 
(In fact, for $\lambda _1,...\lambda _n \in D$, we have
$\| \pi (\lambda _1,...,\lambda _n)f \| = \| \sum _0
^{\infty } a_{k_1...k_n}(f) \lambda _1 ^{k_1}... \lambda _n ^{k_n} \| \le 
\sum _0 ^{\infty } \| a_{k_1...k_n}(f) \| \cdot
| \lambda _1 |^{k_1}...| \lambda _n
| ^{k_n} < \sum _0 ^{\infty } \| a_{k_1...k_n}(f) \|=
\| f \| $, unless all $a_{k_1...k_n}$, for which $k_1+...+k_n>0$, vanish.)
For $B={\Bbb C}$,
we obtain a non-trivial holomorphic family of multiplicative linear
functionals, of constant norm $1$.

\vskip.1cm

Analytic families in Banach spaces, with constant norm (which do not exist in
the scalar valued case, unless they are constant), had been studied,
e.g., in [TW], [GV1], [GV2], [Glo1], [Glo2]. Now we see that their values 
may consist of very natural objects, like surjective homomorphisms (even
multiplicative linear functionals), to which our results are applicable.

\vskip.1cm

{\bf{Remark 3.}}
For Theorems 5 and 6 the spectral hypothesis at the single point $0$
does not extend to any larger set,
like it does for Theorems 3 and
4, cf. Remarks 2 and 4. If in Example 1 we let $B:= {\Bbb{C}}$ and $n:=1$,
then for $\pi (\lambda ) z( \cdot ) := z( \lambda )$ we have, for $z( \lambda
):= \lambda $ and $\lambda =0$, that the spectrum of 
$\pi (\lambda ) z ( \cdot )$ 
for $\lambda =0$ is $\{ 0 \} $, but for any $\lambda \in D 
\setminus \{ 0 \} $ (for Theorem 6 for any $\lambda \in ( -1, 1) \setminus \{
0 \} $), we have that the spectrum of 
$\pi (\lambda ) z( \cdot ) = z (\lambda )  = \lambda \in {\Bbb{C}}$ is
$\{ \lambda \} \ne \{ 0 \} $. See also 
the explanation in Remark 4, after the proof of Theorem 6.

\vskip.1cm

Observe that in Example 1 we could not control the spectra of the elements in
the kernels of $\pi (\lambda _1, \dots , \lambda _n )$. Next we show a
modification of this example, where $\pi (\lambda _1, \dots , \lambda _n )$ is
not constant, but, for the spectra of the elements of Ker\,$\pi (\lambda _1, 
\dots , $
\newline
$ \lambda _n )$, we have the strongest ``smallness'' property from all
our theorems: namely, that all these spectra are $\{ 0 \} $. In this
modification, we may even have 

\newpage

commutative Banach algebras.

\vskip.1cm

{\bf{Example 2.}}
Let us specify $B$ in the above example, as follows. As a Banach space, $B=
C([0,1])$ (of course with complex valued functions). We define the 
\newline
$^*$-operation by conjugation of the functions. Then $B$ will become a Banach 
algebra (without unit), if multiplication is defined as 
the convolution (and also a $^*$-algebra, with
an isometric involution). If $f \in B$, then for its
$n$th convolution power $f^{(*n)}$ we have $\| f^{(*n)} \| 
\le \| f \| ^n /n!$,
hence the spectral radius of any element is $0$, and the spectrum of any
element is $\{ 0 \} $. This readily implies that $B$ has no unit.
(In fact, we could take for $B$ any commutative
radical Banach algebra, of course, without unit.) 
We define $A$, and  $\pi (\lambda _1, \dots , \lambda
_n )$, as in Example 1. 

Then, also for $A$, we have that any element has
spectral radius $0$, and hence, has 
spectrum $ \{ 0 \} $. In fact, for monomials
this is evident. Then, taking into account that, for commutative Banach
algebras, the spectral radius is subadditive, we have that the spectrum of any
polynomial is $\{ 0 \} $, as well. Last, we write any element $f \in A$, 
i.e., an absolutely convergent power series on ${\overline{D}}^n$, 
with coefficients from $B$, as  
$$
f = \sum^\infty_0 \! a_{k_1 \dots k_n} \, z_1^{k_1} \dots z_n^{k_n} =
\sum_{\max k_i \le N} \! a_{k_1 \dots k_n} \, z_1^{k_1} \dots z_n^{k_n} +
\sum_{\max k_i > N} \! a_{k_1 \dots k_n} \, z_1^{k_1} \dots z_n^{k_n}\,.
$$
Here, on the right hand side, the first summand has spectral radius $0$, and
the second summand has spectral radius at most the sum of the respective
norms, which is arbitrarily small, if $N$ is sufficiently large. Hence, the
spectral radius of $f$ is $0$, as asserted. 
Then, of course, for any $(\lambda _1, \dots , \lambda
_n) \in D^n$, we have that the spectra of all elements of Ker\,$\pi (\lambda
_1, \dots , \lambda _n ) \,\,(\subset A)$ 
are $\{ 0 \} $ as well. The only problem is that  
both $B$ and $A$ are not unital.

Therefore, let us consider their unitizations $B \oplus {\Bbb C}$, and  
$A \oplus {\Bbb C}$, with the $l^1$-sum norms, which are Banach algebra norms. 
Then, each $\pi (\lambda _1, \dots , \lambda _n )$ has a
unique extension ${\tilde{\pi}} (\lambda _1, \dots , \lambda _n) :
B \oplus {\Bbb C} \to A \oplus {\Bbb C}$, that is unital: we have 
${\tilde{\pi}} (\lambda _1, \dots , \lambda _n) (f \oplus c):=
\pi (\lambda _1, \dots , \lambda _n )(f) \oplus c$, and 
this extension is a surjective algebra homomorphism as well. Let $f \oplus c 
\in $ Ker\,$
{\tilde{\pi}} (\lambda _1, \dots , \lambda _n)$. Then $0 \oplus 0=
{\tilde{\pi}} (\lambda _1, \dots , \lambda _n)(f \oplus c)=
\pi (\lambda _1, \dots , \lambda _n )(f) \oplus c$, hence $0=c$, and
$f \oplus c = f \oplus 0$. Therefore the spectral radius of $f \oplus c$ is
$0$, and $\sigma (f \oplus c) = \{ 0 \} $.

\vskip.1cm

If we take $n=1$ and $B:={\Bbb C}$ (which is unital), then we have that $A$
from Example 1 is the Banach algebra of absolutely convergent power series on
${\overline{D}}$ (which is also unital, and is not of the above form $B \oplus
{\Bbb{C}}$). Here we have the opposite situation: for $f \in A$ we have 
$\sigma (f)=f({\overline{D}})$, that may disconnect ${\Bbb C}$, even if $f$ is
in the kernel of $\pi (\lambda )$, i.e., if $f(
\lambda )=0$, for some $\lambda \in D$.
We give a concrete example. Let us consider the open $\varepsilon
_0$-neighbourhood $C_{\varepsilon _0}$ 
of the curve $C:=[0,1] \cup \{ z \in {\Bbb{C}} \,\,:\,
|z|=1, $ Im\,$z \ge 0 \} $, where $\varepsilon _0 >0$ is sufficiently small. By
the Riemann mapping theorem there is a bijective analytic map $f: D \to 
C_{\varepsilon _0}$, that
can be extended to a homeomorphism of ${\overline{D}}$ onto 
${\overline{C_{\varepsilon _0}}}$. 
Then $f^2(D)$ will be a small neighbourhood of 
$C^2=[0,1] \cup \{ z \in {\Bbb{C}} \,:\, |z|=1 \} $ (in fact, non-constant
analytic functions are open maps). Hence it contains $0$, and
disconnects ${\Bbb{C}}$. The only problem is that here possibly the power
series of $f$ is not
absolutely convergent. Therefore, let us restrict $f$ to $(1 - \varepsilon
)D$, for $\varepsilon >0$ sufficiently small, whose power series expansion 

\newpage

is
absolutely convergent on $(1 - \varepsilon ){\overline{D}}$.
Then still $f^2\left( (1-
\varepsilon )D \right) $ contains $0$, and $f^2\left( (1-
\varepsilon ){\overline{D}} \right) $ disconnects ${\Bbb{C}}$.
Thus, even the
weakest spectral property from our theorems is not satisfied for this case.

\vskip.1cm

{\bf{Example 3.}}
Unfortunately, in the algebra $A \oplus {\Bbb C}$ from the first part of
Example 2, each element has a one-point
spectrum ($\sigma (f \oplus c) = \{ c \} $), hence this algebra 
has no non-trivial
idempotents. Next we show that all our hypotheses can be fulfilled
simultaneously. Let $A_1$ be any Banach algebra (with continuous involution), 
with a non-central
idempotent $e_1$ (e.g., all complex
$n \times n$ matrices, for $n \ge 2$). Let $x_1
\in A_1$, and $e_1x_1 \ne x_1e_1$. Then, in $A_1$, we have a non-trivial
analytic family of idempotents. Namely, for $\lambda \in {\Bbb C}$, close to
$0$, we have $e_1e^{\lambda x_1} \ne e^{\lambda x_1}e_1$. Therefore, the
element obtained from $e_1$ by a similarity via an exponential function, 
namely $e^{- \lambda x_1}e_1
e^{\lambda x_1}$, is a non-trivial analytic family of idempotents, even defined
for all $\lambda \in {\Bbb C}$. (More generally, 
compare also the theorem of Labrousse, [La], from the title of his paper.)
Now 
we multiply both $B \oplus {\Bbb C}$, and  
$A \oplus {\Bbb C}$, by $A_1$ (taking on the Cartesian
products the $l^{\infty }$-sum
norms).
Then extend ${\tilde{\pi}} (\lambda _1, \dots , \lambda _n)$, as identity on 
$A_1$. Thus we obtain a non-trivial analytic family of surjective 
homomorphisms 
$(A \oplus {\Bbb C}) \times A_1 \to (B \oplus {\Bbb C}) \times A_1$.
Then the elements of the kernels have a component in $A_1$
which is equal to $0$, hence the spectra of these elements are
$\{ 0 \} $. Simultaneously we have a non-trivial analytic family of idempotents
in $(A \oplus {\Bbb C}) \times A_1$; their components in $(A \oplus {\Bbb C})$
are chosen as the identity, or zero, of $A \oplus {\Bbb C}$, and their
components in $A_1$ are the above $e^{- \lambda x_1}e_1e^{\lambda x_1}$.


\heading
4. Proofs of the theorems
\endheading

\demo{Proof of Theorem 1}
We follow \cite{AMMZ, p.~25}, which in turn followed the proof 
of\/ $\pi E(A) = E(A / \text{\rm rad}\, A)$ for 
$\pi : A \to A / \text{\rm rad}\, A$ the canonical map, from \cite{Ri, p.~58} 
and \cite{Ka, p.~124}.

The set of surjective bounded linear maps between Banach spaces is open in 
the corresponding operator space \cite{Gle, Proposition~1.5}, so we may 
suppose that $\pi(\lambda)$ is surjective for each $\lambda \in U$ 
($U$ can be decreased).

By a theorem of Harte, cf.\ \cite{A91, p. 47, Theorem 3.3.8}, 
for any $x \in A$ we have
$$
\sigma(\pi(0)x)
\subset \cap \{\sigma(x + y) : y \in \text{\rm Ker}\, \pi(0)\} \subset
$$
$$
\sigma(\pi(0)x) \cup\bigl( \cup 
\{G : G \text{ is a bounded connected component 
of } \Bbb C \setminus \sigma(\pi(0) x)\}\bigr).
$$
(The last set is also called the {\it polynomially convex hull} 
of\/ $\sigma(\pi(0)x)$.)

Like in the proof of \cite{AMMZ, Theorem~3.1}, we may assume that 
$\sigma(q(0)) = \{0, 1\}$ (else $p \equiv 0$ or $p \equiv 1$ could be chosen).
Applying Harte's theorem, we see that
$$
\cap \{\sigma(c) : c \in \pi(0)^{-1} q(0)\} = \{0, 1\}.
$$
Choose $c \in \pi(0)^{-1} q(0)$ such that $\sigma(c) \not\ni 1/2$.

By \cite{Gle, Lemmas 1.7, 1.10}, there is an open set $V \subset \Bbb C$, 
such that $0 \in V \subset U$, and an analytic map $a : V \to A$ such that 
$a(0) = c$ and
$$
\pi(\lambda) a(\lambda) = q(\lambda) \ \text{ for each } \ \lambda \in V.
$$

\newpage

Hence $\sigma \left( a(0) \right) \not\ni 1/2$.

By upper semicontinuity of the spectrum [A79, p. 6, Th\'eor\`eme 1.1.3] or
\cite{A91, p.~50, Theorem 3.4.2}, 
we may assume, by decreasing~$V$, that
$$
\sigma(a(\lambda)) \not\ni 1/2, \ \text{ for each }\ \lambda \in V.
$$
Then, by the spectral mapping theorem, we have
$$
\sigma(a(\lambda)^2 - a(\lambda)) \not\ni - 1/4, \ \text{ for each }\ 
\lambda \in V.
$$
(For this observe that the only solution of the equation 
$\lambda^2 - \lambda = -1/4$ in $\Bbb C$ is $\lambda = 1/2$.)

Hereafter only small modifications are necessary compared to \cite{Ri} and 
\cite{Ka}.
We have
$$
-a(\lambda)^2 + a(\lambda) =: r(\lambda) \in \text{\rm Ker}\, \pi(\lambda),
$$
with $r : V \to A$ analytic, such that $\sigma(r(\lambda)) \not\ni 1/4$, 
for each $\lambda \in V$.
We want to obtain a solution $z(\lambda) \in \text{\rm Ker}\, \pi(\lambda)$, 
with $z : V \to A$ analytic, of the equation
$$
(a(\lambda) + z(\lambda))^2 = a(\lambda) + z(\lambda),
\tag{1}
$$
that commutes with $a(\lambda)$.

Due to this commutativity, (1) reduces to
$$
z(\lambda)^2 + (2 a(\lambda) - 1) z(\lambda) - r(\lambda) = 0.
\tag{2}
$$
Since $1/2 \notin \sigma(a(\lambda))$, for each $\lambda \in V$, we may 
consider
$$
x(\lambda) := z(\lambda) (2 a(\lambda) - 1)^{-1},
\tag{3}
$$
which is analytic on~$V$.
Observing that $(2a(\lambda) - 1)^2 = 1 - 4r(\lambda)$, we rewrite (2) as an 
equation for $x(\lambda)$, obtaining
$$
x(\lambda)^2 + x(\lambda) - r(\lambda) (1 - 4r(\lambda))^{-1} = 0
\tag{4}
$$
(recall that $\sigma(r(\lambda)) \not\ni 1/4$, for each $\lambda \in V$).
We write
$$
-r(\lambda) (1 - 4r(\lambda))^{-1} =: r_0 (\lambda) \in \text{\rm Ker}\, 
\pi(\lambda),
$$
where $r_0 : V \to A$ is analytic.
Observing that the map
$$
\Bbb C \setminus \{1/4\} \ni \lambda \mapsto - \lambda(1 - 4\lambda)^{-1}
$$
does not contain $1/4$ in its image, the spectral mapping theorem gives
$$
\sigma(r_0(\lambda)) \not\ni 1/4, \ \text{ for each } \ \lambda \in V.
$$

\newpage

We rewrite (4) as
$$
x(\lambda)^2 + x(\lambda) + r_0(\lambda) = 0.
\tag{5}
$$
A formal solution of (5) would be
$$
x(\lambda) = - \frac12 \pm \frac12 \sqrt{1 - 4r_0(\lambda)}.
\tag{6}
$$
Observe that $\sigma(1 - 4r_0(0)) \not\ni 0$.
Moreover, by the hypothesis of the theorem, $\sigma(1 - 4r_0(0))$ does not 
disconnect $\Bbb C$, since $r_0(0) \in \text{\rm Ker}\, \pi(0)$.
Therefore there exists a simple polygonal arc $P$ connecting $0$ with 
infinity, ending with a half-line, and avoiding $\sigma(1 - 4r_0(0))$.
Hence, the distance of $P$ and $\sigma(1 - 4r_0(0))$ is positive. We write
$\varepsilon \in (0, \infty )$ for the third of this distance.

By upper semicontinuity of the spectrum [A79, p. 6, Th\'eor\`eme 1.1.3] or
[A91, p. 9, Theorem 3.4.2], for all $\lambda$ in some open set
$V'$, 
containing~$0$, and contained in $V$, we have that
$$
\cases
\aligned
\sigma(1 - 4r_0(\lambda))
\text{ lies in }
&\text{the} 
{\text{ closed $ \varepsilon $-neighbourhood of }}
\sigma(1 - 4r_0(0)), \\
{\text{thus, in particular, }}
&\text{it} 
{\text{ avoids the closed }} 
\varepsilon {\text{-neighbourhood of }} P.
\endaligned
\endcases
\tag{7}
$$
Denoting by $\rho ( \cdot )$ the spectral radius, we have, by (7),
for all  $\lambda \in V'$, that
$$
\rho \left( 1 - 4r_0(\lambda) \right)  \le
\rho \left( 1 - 4r_0(0) \right) + \varepsilon .
\tag{8}
$$
We may assume that the above open set $V'$ equals $V$ ($V$ may be decreased).

The function $\lambda \mapsto\! \lambda^{1/2}$ has two analytic branches on 
$\Bbb C \setminus P$.
Therefore $\sqrt{1\! - \! 4r_0(\lambda)}$ can be defined by holomorphic 
calculus as
$$
\frac1{2\pi i} \int_\Gamma z^{1/2} \bigl[ z - (1 - 4r_0(\lambda))\bigr] 
^{-1} dz.
\tag{9}
$$
Here $z^{1/2}$ is any of the above two analytic branches.
Moreover, $\Gamma \subset \Bbb C \setminus P$ is a (closed)
Jordan polygon, which is 
the union of two closed arcs $\Gamma _1 $ and $\Gamma_2$, 
with disjoint relative interiors. Here $\Gamma _1$ 
lies in the closed $\varepsilon $-neighbourhood of the simple
polygonal arc~$P$, with both its endpoints having a distance at least 
$\rho \left( 1 - 4 r_0(0) \right) +2 \varepsilon $ from $0$, 
and $\Gamma _2$, having the same endpoints,
has a distance from $0$ at least $\rho \left( 1 - 4r_0(0) \right) +2
\varepsilon $.
(So, $\Gamma _1$ has points close to $0$, as well as points far from $0$, but,
anyway, it lies close to $P$. At the same time, 
$\Gamma _2$ has only points that are far from $0$. A proper
choice for $\Gamma $,
for the case when $P$ is the non-negative real axis, is the following. 
$\Gamma _1$
is the polygonal arc $(1/\varepsilon , - \varepsilon) (0, - \varepsilon) (-
\varepsilon , 0) (0,  \varepsilon) (1/\varepsilon ,  \varepsilon)$ and 
$\Gamma _2 $ is the polygonal arc $(1/\varepsilon ,  \varepsilon) 
(1/\varepsilon , 1/\varepsilon ) (-1/\varepsilon , 1/\varepsilon )
(-1/\varepsilon , -1/\varepsilon ) (1/\varepsilon , -1/\varepsilon )
(1/\varepsilon , -\varepsilon )$. Cf. Fig. 1 in separate file, $\Gamma _1$ in
continuous lines, $\Gamma _2$ in broken lines.)
By (7), (8) and the construction of $\Gamma $, we have that
$$
\Gamma \text{ encloses } \sigma(1 - 4r_0(\lambda)), \text{ for all } \lambda 
\in V.
\tag{10}
$$
Therefore, the definition given in formula (9) is correct.
Then the function 

\newpage

$\sqrt{1 - 4 r_0(\lambda)}$, defined by (9), is analytic 
on~$V$, hence, by (6), $x : V \to A$ is analytic as well, and solves 
equation~(5).
Since $r_0(\lambda)$ is a rational function of\/ $a(\lambda)$, they 
commute, 
hence, by (6) and (9), $x(\lambda)$ and $a(\lambda)$ commute.
Next, by (3), 
$z(\lambda) = (2a(\lambda) - 1) x(\lambda)$ and $a(\lambda)$
commute as well, as we wanted to show.
Since (5) is satisfied, each of (4), (2) and (1) is satisfied, as well.
Lastly,
$$
\pi(\lambda)(a(\lambda) + z(\lambda)) = \pi(\lambda) (a(\lambda) 
+ (2a(\lambda) - 1) x(\lambda)) = q(\lambda) + (2q(\lambda) - 1) 
\pi(\lambda)(x(\lambda)),
$$
and, by (6) and (9), with $\Gamma$ satisfying (10), we have
$$
\aligned
\pi(\lambda)(x(\lambda))
&= -\frac12 \pm \frac12 \cdot \frac1{2\pi i} \int_\Gamma z^{1/2} 
\big[  z - \pi(\lambda) \left( 1 - 4 r_0(\lambda ) \right) \big] ^{-1} 
dz\\
&= -\frac12 \pm \frac12 \cdot \frac1{2\pi i} \int_\Gamma z^{1/2} 
(z - 1)^{-1} dz
\endaligned
\tag{11}
$$
(since $\pi(\lambda) r_0(\lambda) = 0$).
Since $\Gamma$, by (10), encloses $\sigma(1 - 4 r_0(\lambda))$, for each 
$\lambda \in V$, therefore it encloses also its subset 
$\sigma\bigl[ \pi(\lambda)(1 - 4r_0(\lambda))\bigr] = \{1\}$.
Then
$$
\frac1{2\pi i} \int_\Gamma z^{1/2} (z - 1)^{-1} dz
$$
is already independent of\/ $\lambda$, depends only on $\Gamma$ and the 
chosen branch of\/ $\lambda \mapsto \lambda^{1/2}$, and has a value $1$ 
or $-1$.
Therefore, by (11),
$$
\pi(\lambda)(x(\lambda)) = - \frac12 \pm \frac12(\pm 1),
$$
hence, with proper choice of the first $\pm$ sign, we have 
$\pi(\lambda) (x(\lambda)) = 0$.
Thus, by (3), we have
$$
\pi(\lambda) (z(\lambda)) = 0,
$$
therefore
$$
\pi(\lambda) (a(\lambda) + z(\lambda)) = q(\lambda),
$$
while, by (1), that we already know to hold, $a(\lambda) + z(\lambda)$ is 
idempotent, moreover $a(\lambda) + z(\lambda)$ is analytic on~$V$.\hfill 
$\blacksquare$
\enddemo


Before the proof of Theorem 2 we recall some elementary concepts from
topology. Total disconnectedness of a topological space was recalled in the
first paragraph of Section 2. 
A topological space is {\it{$0$-dimensional}} if it is
non-empty, and has an
open base consisting of open-and-closed sets. In particular, the totally
disconnected spectra, mentioned in the first paragraph of Section 2, are
$0$-dimensional (cf. {\bf{P1}} below).
For these spaces the following
two properties are well known.

{\bf{P1}}. Total disconnectedness and $0$-dimensionality are equivalent for
{\it{non-empty compact Hausdorff spaces}} [E, p. 362, Theorem 6.2.10, and 
p. 388, Theorem 7.1.12], and, for the non-empty compact 
metric case, [Ku, p. 189, Section 47, IX, 2nd paragraph] 
(but they are not equivalent for general spaces, 
cf. [Ku, p. 152, Section 46, VI, Remark (i)]). (Observe that
both [Ku, p. 151, Section 46, VI, Definition]

\newpage

and [E, pp. 356, 360, 369] 
call total disconnectedness ``hereditary disconnectedness'', which is a rather
obsolete terminology. What [Ku, p. 151, Section 46, VI, Definition] 
calls total disconnectedness, is another property, not used in our paper.)
We note that we
use these properties only for spectra of elements in Banach algebras, which
are non-empty and compact.

{\bf{P2}}. 
The complement of a $0$-dimensional subset of the plane is connected 
[Ku, p. 188, Section 47, VIII, Theorem 1, and p. 466, Section 
59, II, Theorem 1].


\demo{Proof of Theorem 2}
Like in the proof of Theorem~1, we may suppose that $\sigma(q(0)) = \{0, 1\}$.
By analyticity, both $\pi$ and $q$ extend to analytic functions 
$\pi : U \to \Cal B(A, B)$ and $q : U \to E(A)$ for some open set $U$ being
the union of small discs about each point of $G$, 
with $G \subset U \subset G + i \Bbb R \subset \Bbb C$, 
hence with $U \cap \Bbb R = G$, and with each connected component of\/ $U$ 
intersecting~$G$.
We have that $\pi(\lambda)$ is a homomorphism $A \to B$.
Both the homomorphism property of\/ $\pi(\lambda)$, and the 
idempotency of\/ $q(\lambda)$ are consequences of the identity theorem for 
analytic functions. 
Like in the proof of Theorem~1, we may suppose that $\pi(\lambda)$ is 
surjective for each $\lambda \in U$.

As in the proof of Theorem~1, there exists an open set $V \subset \Bbb C$ 
such that $0 \in V \subset U$, and an analytic map $a : V \to A$, such that
$$
\pi(\lambda) a(\lambda) = q(\lambda) \ \text{ for each } \ \lambda \in V.
$$
Let $\lambda \in V \cap \Bbb R$.
Since $\pi(\lambda)$ is a $^*$-homomorphism, therefore
$$
\pi(\lambda) (a(\lambda)^*) = (\pi(\lambda) a(\lambda))^* = q(\lambda)^* 
= q(\lambda) = \pi (\lambda) a(\lambda).
$$
This implies
$$
q(\lambda) = \pi(\lambda) \frac{a(\lambda) + a(\lambda)^*}{2} \ 
\text{ for each } \ \lambda \in V \cap \Bbb R.
\tag{12}
$$
Here $a(\lambda)$ is analytic, $a(\lambda)^*$ is conjugate analytic, 
so their restrictions to $V \cap \Bbb R$ are real analytic.
(Observe, that in the interior of a circle of convergence, we have that
$$
a(\lambda) = \sum^\infty_0 a_n(\lambda - \lambda_0)^n.
$$
By continuity of the involution on $A$, this implies that
$$
a(\lambda)^* = \sum^\infty_0 a_n^* (\ol \lambda - \ol \lambda_0)^n.
$$
So, in fact, $a(\lambda)^*$ is conjugate analytic on each open set 
of~$\Bbb C$, where $a(\lambda)$ is analytic.)

Thus we have a real analytic, self-adjoint valued function
$$
a_0(\lambda) := \frac{a(\lambda) + a(\lambda)^*}{2}
\tag{13}
$$
on $V \cap \Bbb R$, lifting $q(\lambda)$,
for each $\lambda \in V \cap \Bbb R$, by~(12).
This function $a_0(\lambda)$ extends 

\newpage

by analyticity to an analytic function 
$a_0(\lambda)$ on a neighbourhood of\/ $V \cap \Bbb R$, contained 

in 
$(V \cap \Bbb R) + i \Bbb R \subset \Bbb C$.
We may suppose that this neighbourhood is $V$ ($V$ may be decreased).
By analyticity of\/ $\pi(\lambda)$, $q(\lambda)$ and $a_0(\lambda)$, and the 
identity theorem for analytic functions, we have, by (12) and (13), for each 
$\lambda \in V$, that
$$
\pi(\lambda) a_0(\lambda) = q(\lambda).
$$
This means that we may assume that our original lifting $a(\lambda)$ is 
self-adjoint valued on $V \cap \Bbb R$.
Then, for each $\lambda \in V \cap \Bbb R$, we have
$$
\sigma(a(\lambda)) = \ol{\sigma(a(\lambda))}.
$$
(Observe that in general we do not have that the spectrum of a self-adjoint
element is real: in the Banach algebra of ${\Bbb C}$-valued absolutely
convergent power series on $\overline{D}$, with norm $\| \sum _{n=0} ^ {\infty
} c_n z^n\| :=\sum _{n=0} ^ {\infty} |c_n|$, and with involution given by
coefficientwise conjugation, the self-adjoint element $f(z) \equiv z$ has
spectrum $\overline{D}$.)

Then, like in the proof of \cite{AMMZ, Theorem 3.1}, we construct the Riesz 
idempotent $p(\lambda) \in A$, and elements $a_0(\lambda), 
a_1(\lambda) \in A$. These functions are analytic on a neighbourhood of\/ $0$; 
we may suppose that this neighbourhood is $V$ ($V$ may be decreased).
Namely,
$$
\align
p(\lambda) &:= \frac1{2\pi i} \int_{\Gamma_1} (z - a(\lambda))^{-1}\, dz,\\
a_0(\lambda) &:= \frac1{2\pi i} \int_{\Gamma_0} 
(1 - z)^{-1} (z - a(\lambda))^{-1}\, dz,\\
a_1(\lambda) &:= \frac1{2\pi i} \int_{\Gamma_1} z^{-1}(z - a(\lambda))^{-1} \, 
dz,
\endalign
$$
where $\Gamma_0$ and $\Gamma_1$ are Jordan polygons, their interior domains 
Int\,$\Gamma _i$ being disjoint,
and together covering $\sigma \left( a(0) \right) $, 
with $0 \in \text{\rm Int}\, \Gamma_0$ 
and $1 \in \text{\rm Int}\, \Gamma_1$. These exist, as asserted in [AMMZ,
p. 25], as soon as we know that $\sigma \left( a(0) \right) $ 
is totally disconnected --- that will be proved in the next paragraph --- 
due to the following facts. 

1) The set $\sigma \left( a(0) \right) $ 
is the union of two of its relatively
open-and-closed subsets, one containing $0$, the other one containing $1$, 
which follows from
the fact that total disconnectedness of $\sigma \left( a(0) \right) $ 
implies that $\sigma \left( a(0) \right) $ is $0$-dimensional, cf. {\bf{P1}},
and hence ${\Bbb C} \setminus \sigma \left( a(0) \right) $ 
is connected, cf. {\bf{P2}}. 

2) Connected open subsets of the plane are connected via polygonal arcs 
[Ku, p. 461, Section 59, I, Theorem 1]. 

Now we show that $\sigma \left( a(0) \right) $ is totally disconnected (as
mentioned in [AMMZ, p. 25] without proof).
In fact, we have $\pi (\lambda ) \left( a( \lambda ) ^2 - 
a ( \lambda ) \right) =q ( \lambda )^2- q ( \lambda ) = 0$, 
for each $\lambda \in V$. Hence, in particular, 
$$
a(0)^2-a(0) \in {\text{Ker\,}}\pi (0).
$$
By our spectral hypothesis, $\sigma \left( a(0)^2-a(0) \right) \subset
{\Bbb{C}}$ is totally disconnected. 

\newpage

If $\sigma \left( a(0) \right) $
$ \subset {\Bbb{C}}$ were not totally disconnected, it would
contain, by definition (recalled in the 
\newline
first paragraph of Section 2),
a connected subset 
$$
S \subset \sigma \left( a(0) \right)
$$ 

consisting of more than one points. Then the image
of $S$, under the continuous map $f:{\Bbb{C}} \to {\Bbb{C}}$, defined by
$f(z):= z^2-z$, i.e., the set $f(S)$, 
is connected as well.
Clearly, the connected set $S \subset {\Bbb{C}}$ 
cannot consist of exactly two points. Therefore, let $s_1,s_2,s_3 \in S$
be distinct. Now observe that the inverse image of any $w
\in {\Bbb{C}}$ under $f$ 
consists of at most two points. Therefore, the subset
$\{ f(s_1), f(s_2), f(s_3) \} $ of $f(S)$ consists of at least two points.
Now, also using the spectral mapping theorem, we have
$$
\{ f(s_1), f(s_2), f(s_3) \} \subset 
f(S) \subset f \left( \sigma \left( a(0) \right) \right) =
\sigma \left( f \left( a(0) \right) \right) =
\sigma \left( a(0)^2-a(0) \right) .
$$
Hence, $f(S)$ is a connected subset of $\sigma \left( a(0)^2-a(0) \right) $,
consisting of more than one points, 
while $\sigma \left( a(0)^2-a(0) \right) $ is totally
disconnected, that is a contradiction.
Hence, $\sigma \left( a(0) \right) $ is totally disconnected, as well, 
as promised to be shown at the beginning of this paragraph.

By upper semicontinuity of the spectrum [A79, p. 6, Th\'eor\`eme 1.1.3] or
\cite{A91, p.\ 50, Theorem 3.4.2}, we may 
suppose that this decomposition of the spectrum by the Jordan polygons $\Gamma
_0$ and $\Gamma _1$ (i.e., being covered by $($Int\,$\Gamma _0) \cup (
$Int\,$\Gamma _1)$) holds for all $a(\lambda)$, 
where $\lambda \in V$ ($V$ may be decreased).
Like in the proof of \cite{AMMZ, Theorem 3.1}, we obtain
$$
a(\lambda) - p(\lambda) = 
(a(\lambda)^2 - a(\lambda))(a_1(\lambda) - a_0(\lambda)) 
\in \text{\rm Ker}\, \pi(\lambda).
$$
(For the reader's convenience, we sketch its proof.
We have $p(\lambda) = a(\lambda) a_1(\lambda)$, and $1 - p(\lambda) 
= (1 - a(\lambda)) a_0(\lambda)$, hence $a(\lambda) - p(\lambda) 
= a(\lambda)(1 - p(\lambda)) - (1 - a(\lambda)) p(\lambda) 
= a(\lambda)(1 - a(\lambda))a_0(\lambda) - (1 - a(\lambda)) a(\lambda) 
a_1(\lambda)$.)
Therefore
$$
\pi(\lambda) p(\lambda) = \pi(\lambda) a(\lambda) = q(\lambda) \ 
\text{ for each }\ \lambda \in V.
$$

It remains to ensure that $p(\lambda) = p(\lambda)^*$ for $\lambda \in H 
:= V \cap \Bbb R$.
However, by $\sigma \left( a(0) \right) = \overline{\sigma \left( a(0) \right)
}$,
we can choose $\Gamma_1$ symmetric with respect to the real axis, which yields 
$p(\lambda) = p(\lambda)^*$, for each $\lambda \in H$.\hfill $\blacksquare$
\enddemo


\demo{Proof of Theorem 3}
The proof is analogous to those of \cite{Ka, pp.\ 125--126, Theorem 31, 
Corollary}. For convenience, we use the notation of~\cite{Ka}.

We only have to prove the statement corresponding to \cite{Ka, p.~125, 
Theorem 31}, since \cite{Ka, p.\ 126, Corollary} is an easy consequence, 
whose proof can be taken over without modification (cf. the last paragraph of
this proof).

Thus we only have to prove the following statement.
Let $u, v: U \to E(B)$ be analytic, with $u(\lambda)$ and $v(\lambda)$ 
orthogonal, for each $\lambda \in U$.
Let $e: U \to E(A)$ 

\newpage

be analytic, with $\pi(\lambda) e(\lambda) 
= u(\lambda)$ for each $\lambda \in U$.
Then there exists $f: U \to E(A)$ analytic, such that
$$
\cases
\aligned
e(\lambda)  
\text{ and } 
&f(\lambda) = f(\lambda)^2 \text{ are orthogonal 
for each }\lambda \in U,\\ 
&\text{and }\pi(\lambda) f(\lambda) = v(\lambda)
\text{ for each }\lambda \in U.
\endaligned 
\endcases
 \tag{14}
$$

By \cite{Le1, Theorem 5.1} there exists $b : U \to A$ analytic, such that
$$
\pi(\lambda) b(\lambda) = v(\lambda) \ \text{ for each }\ \lambda \in U.
$$
Then, following \cite{Ka}, we define
$$
a(\lambda) := (1 - e(\lambda)) b(\lambda)(1 - e(\lambda)).
$$
Then
$$
\pi(\lambda) a(\lambda) = v(\lambda),
$$
and $e(\lambda) a(\lambda) = a(\lambda) e(\lambda) = 0$, and
$$
z(\lambda) := a(\lambda)^2 - a(\lambda) \in \text{\rm Ker}\, \pi(\lambda),
$$
with $e(\lambda)z(\lambda) = z(\lambda) e(\lambda) = 0$.

Then, by hypothesis, 
$$
\sigma \left( z( \lambda ) \right) = \{ 0 \}, 
\tag{15}
$$
and therefore
$\left( 2 a(\lambda) - 1 \right) ^2 = 1 + 4z( \lambda ) $ is invertible.
Also, by the spectral mapping theorem, we have by (15) that
$$
\sigma \left( a(\lambda) \right) \subset \{ 0,1 \}.
\tag{16}
$$

Then, following \cite{Ka}, we solve the equation
$$
w(\lambda)^2 + w(\lambda) + z(\lambda)(2a(\lambda) - 1)^{-2} = 0,
\tag{17}
$$
where $w : U \to A$ is a function.
One of the formal solutions of (17) is
$$
w(\lambda) = -\frac12 + \frac12 \sqrt{1 - 4z(\lambda)(2 a(\lambda) - 1)^{-2}}.
$$
Here 
$$
4z(\lambda)(2a(\lambda) - 1)^{-2} \in \text{\rm Ker}\, \pi(\lambda),
$$
hence 
$$ \sigma \left( 4z(\lambda)(2a(\lambda) - 1)^{-2} \right) 
=\{0\}.
\tag{18}
$$

\newpage

Therefore, by functional calculus, $w(\lambda)$ can be defined as
$$
w(\lambda) := - \frac12 + \frac12 \frac1{2\pi i} \int_\Gamma (1 - \zeta)^{1/2}
\bigl[ \zeta - 4z(\lambda)(2 a(\lambda) - 1)^{-2}\bigr] ^{-1}\, d\zeta,
\tag{19}
$$
where $\Gamma$ is a small circle with centre at $0$, and we take that branch 
of\/ $(1 - \zeta)^{1/2}$, that takes some positive values on~$\Gamma$.
Then $w : U \to A$ is analytic, and solves 
equation~(17), and $w(\lambda)$ 
commutes with $a(\lambda)$ and $e(\lambda)$.
Moreover, $w(\lambda) \in \text{\rm Ker}\, \pi(\lambda)$, since
$$
\align
\pi(\lambda) w(\lambda)
&= - \frac12 + \frac12 \cdot \frac1{2\pi i} \int_\Gamma (1 - \zeta)^{1/2} 
\bigl[\zeta - \pi(\lambda)\left( 4z(\lambda)(2a(\lambda) - 1)^{-2} \right)
\bigr] ^{-1}\, 
d\zeta\\
&= - \frac12 + \frac12 \cdot \frac1{2\pi i} 
\int_\Gamma (1 - \zeta)^{1/2}(\zeta - 0)^{-1}\, d\zeta 
= -\frac12 + \frac12 \cdot 1 = 0.
\endalign
$$

Let
$$
x(\lambda) := (1 - e(\lambda)) w(\lambda) 
= w(\lambda)(1 - e(\lambda)) \in \text{\rm Ker}\, \pi(\lambda),
$$
with $x : U \to A$ analytic.
Then $e(\lambda) x (\lambda) = x(\lambda) e(\lambda) = 0$.
Let $r : U \to A$ be the analytic function, defined by
$$
r(\lambda) := x(\lambda) (2a(\lambda) - 1) 
= (2a(\lambda) - 1) x(\lambda) \in \text{\rm Ker}\, \pi(\lambda).
$$
Note that $r(\lambda )$ commutes with $a(\lambda )$. (Actually, all elements
$a(\lambda ), e( \lambda ), z(\lambda ), w(\lambda )$, 
\newline
$x(\lambda ), r(\lambda
)$ lie in a commutative subalgebra of $A$, depending on $\lambda $,
as can be seen from their definitions, step by step. Despite the delicate fact
that $a ( \lambda )$'s, for different $\lambda $'s, may not be commuting!)

Let $f : U \to A$ be the analytic function, defined by
$$
f(\lambda) := a(\lambda) + r(\lambda).
$$
Then, like in \cite{Ka} (by essentially straightforward calculations, 
using commutativity of $a(\lambda ), r(\lambda )$, and (17)), we obtain
$$
f(\lambda)^2 = f(\lambda), \text{ and }\ e(\lambda) f(\lambda) 
= f(\lambda) e(\lambda) = 0.
\tag{20}
$$
Finally,
$$
\pi(\lambda) f(\lambda) = \pi(\lambda) a(\lambda) = v(\lambda).
\tag{21}
$$

\newpage

This ends the proof of (14).

Now the proof of \cite{Ka, p.~126, Corollary} finishes the proof of our 
theorem.
(For the reader's convenience, we recall that, to define 
$p_k$ if\/ $p_1, \dots p_{k - 1}$ already were defined, (14) was there 
applied for $u := q_1 + \dots + q_{k - 1}$, and $v := q_k$, and 
$e := p_1 + \dots + p_{k - 1}$; here, for $k=1$, an empty sum means $0$.)
\hfill $\blacksquare$
\enddemo


\demo{Proof of Theorem 4}
Again we proceed like \cite{Ka, pp.\ 125--126, Theorem 31, Corollary}.
Once more, we only have to prove the statement corresponding to 
\cite{Ka, p.~125, Theorem 31}, since \cite{Ka, p.\ 126, Corollary} is an 
easy consequence, whose proof can be taken over without modification (cf. the
last paragraph of this proof).

Thus we only have to prove the following statement.
Let $u, v : G \to S(B)$ be real analytic, with $u(\lambda)$ and 
$v(\lambda)$ orthogonal, for each $\lambda \in G$.
Let $e : G \to S(A)$ be real analytic, such that $\pi(\lambda) e(\lambda) 
= u(\lambda)$ for each $\lambda \in G$.
Then there exists $f : G \to S(A)$ real analytic, such that
$$
\cases
\aligned
e(\lambda) \text{ and } f(\lambda)
&= f(\lambda)^2 = f(\lambda)^* \text{ are orthogonal for each } 
\lambda \in G,\\
&\text{and } \pi(\lambda) f(\lambda) = v(\lambda) \text{ for each } 
\lambda \in G.
\endaligned
\endcases
\tag{22}
$$

Like in the proof of Theorem~2, by analyticity, $\pi, u, v$ and $e$ extend 
to analytic functions $\pi : U \to \Cal B(A, B)$, and $u, v : U \to E(B)$, 
and $e : U \to E(A)$ for some open set $U \subset \Bbb C$ such that 
$G \subset U \subset G + i \Bbb R \subset \Bbb C$, thus $U \cap \Bbb R = G$, 
with each connected 
component of\/ $U$ intersecting $G$, and with the following properties.
We have that $\pi(\lambda)$, for each $\lambda \in U$, is a surjective 
homomorphism $A \to B$, and, for each $\lambda \in G$, each 
$x \in \text{\rm Ker}\, \pi(\lambda)$ satisfies $\sigma(x) = \{0\}$.
Further, $u(\lambda)$ and $v(\lambda)$ are orthogonal for each 
$\lambda \in U$, and $\pi(\lambda)e(\lambda) = u(\lambda)$ for each 
$\lambda \in U$ (again, by using the identity theorem for analytic functions).
Also at each further step, when the open set $U \subset \Bbb C$, 
containing $G$, is decreased, we always suppose that each connected 
component of~$U$ intersects~$G$.

Then we are in the situation of Theorem~3, except that we do not have 
$\sigma(x) = \{0\}$ for each $x \in \text{\rm Ker}\, \pi(\lambda)$, and for 
each $\lambda \in U$. (However, see Remark 4 after the proof of Theorem 6!)
We have
used this hypothesis, for $x=z( \lambda )$, i.e., $\sigma \left( z( \lambda
) \right) = \{0\}$, 
cf. (15), to prove that 
$\sigma(a(\lambda))\subset \{0, 1\}$ --- and hence $(2a(\lambda) - 1)^{-2}$
exists --- and
$\sigma\bigl(4z(\lambda)(2a(\lambda) - 1)^{-2}\bigr) = \{0\}$, cf. (16), 
(18).
By upper semicontinuity of the spectrum [A79, p. 6, Th\'eor\`eme 1.1.3] or
[A91, p. 50, Theorem 3.4.2], and possibly decreasing 
$U$ at the respective steps of the proof, we may suppose the following three
facts (23), (24), (25) (where always $\lambda \in U$).

After possibly decreasing $U$, 
we may suppose, for some sufficiently small
$\varepsilon >0$, that for all $\lambda \in U$, we have
$$
\sigma \left( z ( \lambda ) \right) \subset \{ z \in {\Bbb{C}} : |z| <
\varepsilon \} . 
\tag{23} 
$$
Namely, by upper semicontinuity of the spectrum, each
point in $G$ has an open neighbourhood, contained in $U$, on which (23)
holds. Then the new (decreased) $U$ will be the union of all these
neighbourhoods.
Therefore, by the spectral mapping theorem, we have, for each $\lambda \in U$,
rather than (16),

\newpage

$$
\sigma(a(\lambda)) \subset 
\{z \in \Bbb C : |z| < 1/3 {\text{ or }} |1 - z| < 1/3\}.
\tag{24}
$$
Then (24) implies that
$(2a(\lambda) - 1)^{-2}$ exists, for each $\lambda \in U$. 
Next, 
possibly decreasing $\varepsilon >0$ and $U$, 
similarly as above, we may suppose, for each $\lambda \in
U$, rather than (18), that
$$
\sigma \bigl(4z(\lambda)(2a(\lambda) - 1)^{-2}\bigr) \subset 
\{z \in \Bbb C : \, |z| < 1/3\}.
\tag{25}
$$
Furthermore, we choose $\Gamma $ in (19) as $\{ z \in {\Bbb{C}} : |z| = 1/2 \}
$.

Then the proof of (14) in the proof of Theorem~3 goes through.
Therefore, there exists $f : U \to E(A)$ analytic, such that 
$e(\lambda)$ and $f(\lambda)$ are orthogonal, for each $\lambda \in U$, 
and $\pi(\lambda) f(\lambda) = v(\lambda)$, for each $\lambda \in U$.
The only property that remains to be ensured is that $f(\lambda)$ is 
self-adjoint, for each $\lambda \in U \cap \Bbb R = G$.
Of course, the proof of Theorem~3 does not give this.
We will have to go over the proof of Theorem~3, and step by step we will 
have to do the respective modifications, that finally will ensure that 
$f(\lambda)$ is self-adjoint, for each $\lambda \in G$.
Like in the proof of Theorem~3, there exists $b : U \to A$ analytic, 
such that $\pi(\lambda) b(\lambda) = v(\lambda)$ for each $\lambda \in U$.
Like in the proof of Theorem~2, we even may assume (by changing $b(\lambda)$) 
that $b(\lambda)$ is self-adjoint, for each $\lambda \in U \cap \Bbb R = G$.
(At the same time possibly $U$ has to be decreased, still each connected
component of\/ $U$ intersecting~$G$.)

Then $a(\lambda)$ is self-adjoint for each $\lambda \in G$, as well as 
$z(\lambda)$. 
Then,
like in (23), (24), (25), we may assume, after possibly decreasing $U$,
for each $\lambda \in U$, the following facts --- 
which are actually of the form
(23), (24), (25), but now for the changed, already self-adjoint
$b(\lambda )$, etc.
For some sufficiently small $\varepsilon >0$, we have
$$
\sigma \left( z ( \lambda ) \right) \subset \{ z \in {\Bbb{C}} : |z| <
\varepsilon \} , 
$$
and hence $\sigma(a(\lambda)) \subset \{z \in \Bbb C : |z| < 1/3$ or 
$|1 - z| < 1/3\}$ --- and hence $(2a(\lambda) - 1)^{-2}$ exists --- and
$\sigma\bigl(4z(\lambda)(2 a (\lambda) - 1)^{-2}\bigr) \subset 
\{z \in \Bbb C : |z| < 1/3\}$.
Then, taking $\Gamma := \{z\in \Bbb C : |z| = 1/2\}$ in (19), we have that 
$w(\lambda)$, defined by (19) and this $\Gamma$, is self-adjoint, for each 
$\lambda \in U \cap \Bbb R = G$.
Observe that the product of two commuting self-adjoint elements is 
self-adjoint.
Therefore, $x(\lambda)$ is self-adjoint, for each $\lambda \in G$, 
and 
then $r(\lambda)$ is self-adjoint, for each $\lambda \in G$, hence, finally, 
$f(\lambda)$ is self-adjoint, for each $\lambda \in G$, as was to be proven.

Now the analogue of the proof of \cite{Ka, p.~126, Corollary} (cf. the last
paragraph of the proof of Theorem 3) finishes the 
proof of our theorem.\hfill $\blacksquare$
\enddemo


\demo{Proof of Theorem 5}

\newpage

We follow the proofs of Theorems 3 and 4. 

As in the proof of Theorem 1, we may assume that $\pi ( \lambda )$ is
surjective for each $\lambda \in U$.

It suffices to prove the analogue of the statement in the third paragraph
of the proof of Theorem 3, where however, rather than $f:U \to E(A)$, we have
only $f:V \to E(A)$, for some open set $V \subset {\Bbb C}$, such that $0 \in
V \subset U$, and (14) holds only for each $\lambda \in V$ (at both places).

As in the proof of Theorem 4, 
like in (23), (24), (25), we may assume, after possibly decreasing $V\,\,(\ni
0)$ at the
respective steps, for each $\lambda \in V$, the following facts.
For some sufficiently small $\varepsilon >0$, we have
$$
\sigma \left( z ( \lambda ) \right) \subset \{ z \in {\Bbb{C}} : |z| <
\varepsilon \} , 
$$
so that $\sigma(a(\lambda)) \subset \{z \in \Bbb C : |z| < 1/3$ or 
$|1 - z| < 1/3\}$ --- hence $(2a(\lambda) - 1)^{-2}$ exists --- and
$\sigma\bigl(4z(\lambda)(2 a (\lambda) - 1)^{-2}\bigr) \subset 
\{z \in \Bbb C : |z| < 1/3\}$.
Then $\Gamma $ can be chosen in (19) as $\{ z \in {\Bbb C} 
: |z| = 1/2 \} $. With these changes, we arrive at formulas (20) and (21),
however, only for $\lambda \in V$.

Last, the proof of [Ka, p. 126, Corollary] (cf. the last paragraph of the
proof of Theorem 3) finishes the proof of our theorem. 
(Observe that if we would have had an infinite sequence $q_1,q_2, \dots $ in
the theorem, then
the neighbourhoods of $0$ in the inductive proof could shrink to
$0$, but in finitely many steps this cannot occur.) \hfill $\blacksquare$
\enddemo


\demo{Proof of Theorem 6}
We follow the proofs of Theorems 4 and 5. 
As in the proof of Theorem 1, we may assume that $\pi ( \lambda )$ is
surjective for each $\lambda \in G$.
It suffices to prove the analogue of the statement in the second paragraph
of the proof of Theorem 4, where however, rather than $f:G \to S(A)$, we have
only $f:H \to S(A)$, for some open set $H \subset {\Bbb R}$, such that $0 \in
H \subset G$, and (22) holds only for each $\lambda \in H$ (at both places).

As in the proofs of Theorems 4 and 5, 
like in (23), (24), (25), we may assume, after possibly decreasing $V\,\,(\ni
0)$ at the
respective steps, for each $\lambda \in V$, the following facts.
For some sufficiently small $\varepsilon >0$, we have
$$
\sigma \left( z ( \lambda ) \right) \subset \{ z \in {\Bbb{C}} : |z| <
\varepsilon \} , 
$$
and hence $\sigma(a(\lambda)) \subset \{z \in \Bbb C : |z| < 1/3$ or 
$|1 - z| < 1/3\}$ --- so that $(2a(\lambda) - 1)^{-2}$ exists --- and
$\sigma\bigl(4z(\lambda)(2 a (\lambda) - 1)^{-2}\bigr) \subset 
\{z \in \Bbb C : |z| < 1/3\}$.
Then $\Gamma $ can be chosen

\newpage

 in (19) as $\{ z \in {\Bbb C} 
: |z| = 1/2 \} $. With these changes, 
also changing $b (\lambda )$ to a self-adjoint element, for $\lambda \in H$,
as in the proofs of Theorems 2 and 4, 
we arrive at (22),
however, only for $\lambda \in H$ (at both places), where $H \subset {\Bbb R}$
is some open set, such that 
$0 \in H \subset G$.

Last, the analogue of the proof of [Ka, p. 126, Corollary] 
(cf. the last paragraph of the proof of Theorem 3) finishes the proof 
of our theorem. \hfill $\blacksquare$
\enddemo


{\bf{Remark 4.}}
Since, in the proof of Theorem 4, we have  
$\sigma \left( z(\lambda ) \right) = \{ 0 \} $ for all $\lambda \in G$
(by $\pi (\lambda ) z(\lambda ) = 0$),
it follows that $\sigma \left( z(\lambda ) \right) = \{ 0 \} $ for all
$\lambda \in U$. Indeed, the logarithm of the spectral radius of $z( \lambda
)$ is a subharmonic function on $U$, by 
Vesentini's theorem (see, e.g., 
[A79, p. 9, Th\'eor\`eme 1.2.1], [A91,
p. 52, Theorem 3.4.7] or [Ra, p. 178, Theorem 6.4.2]). So, the above claim
follows from [AG, p. 125, Corollary 5.1.5 (i)]
(the definition of a {\it{polar set}} mentioned there cf. with [HK, p. 212]), 
or, alternatively, from [A79, p. 11, Corollaire 1.2.1] with [A91, p. 180,
Corollary A.1.27],
[Ra, p. 42,
Exercise 2.5.1], or their strengthenings [A91, p. 180, Theorem A.1.28], 
[Ra, p. 57, Exercise 3.2.1], 
and [HK, p. 225, Theorem 5.13] (the last one
implies that the ``non-degenerate straight line 
segments'' in the 
last sentence of this paragraph
can be replaced by compact sets of any positive Hausdorff dimension),
applied to each connected component of $U$. 
The facts cited in the preceding sentence are particular cases of
H. Cartan's theorem,
stating, for a connected open set $\emptyset \ne U \subset {\Bbb{C}}$, that 

1) if a subharmonic function $f:U \to
[- \infty , \infty )$ is not identically $- \infty $ on $U$, then 
$f^{-1}(- \infty )$ is a $G_{\delta }$ set of capacity zero, and, conversely,

2) if $E \subset U$ is of zero capacity, then there exists 
a subharmonic function $f:U \to
[- \infty , \infty )$, not identically $- \infty $ on $U$, such that 
$E \subset f^{-1}(- \infty )$,
\newline
cf., e.g., [HK, p. 274, Theorem 5.32], 
[A79, p. 173, Th\'eor\`eme 14], and [AZr, proof of Th\'eor\`eme 2.2].
(Concerning implication 1), see also [A91, p. 180, Theorem A.1.29] and
[Ra, p. 65, Theorem 3.5.1], and its weaker form
[Ra, p. 41, Corollary 2.5.3].)
H. Cartan's theorem shows that the strongest spectral assumption 
on Ker\,$\pi ( \lambda )$, in Theorems 3 and 4, 
can be postulated merely for $\lambda $'s, e.g., 
on non-degenerate straight line segments in each
connected component of $U$, or $G$.

Actually, both for Theorem 3 and Theorem 4, the strongest spectral assumption 
on Ker\,$\pi ( \lambda )$ can be postulated
on any subsets $E$ of each connected component of 
$U$, or $G$, respectively, of positive outer 
capacity, but it is not sufficient to postulate them for subsets $E$ of zero
capacity, if we only use
H. Cartan's theorem for general subharmonic functions, as we have seen in this
remark above, and especially in our earlier Remark 3.

\newpage

Cf. also Deny's theorem, [HK, p. 274] or [Ra, p. 65],
which
is a sharper form of H. Cartan's theorem above, and gives a complete
characterization. Namely, if $E \subset {\Bbb{C}}$, and
$\emptyset \ne U \subset {\Bbb{C}}$ 
is any connected open set containing $E$, then the following 
are equivalent:

1) there exists a subharmonic
function $f:U \to [- \infty , \infty )$, not identically $- \infty $,
such that $f^{-1}( -\infty )=E$; 

2) $E$ is a $G_{\delta }$-set of capacity zero. 
\newline
This is also a generalization
of Evans' theorem [A91, p. 179, Theorem A.1.24]
--- that is the particular case of Deny's theorem, for
$E$ compact --- used in the proofs of [AZr, Th\'eor\`eme 2.2] and
[AMZ, p. 520, Lemma]. Let us
also mention that we observed later
that [AMZ, p. 520, Lemma] had been given already earlier, 
even in a more general form, as
[AZr, Th\'eor\`eme 2.2], by the first mentioned author of [AMZ] and
Zra\"\i bi.


\heading
5. Some problems
\endheading

1) Does the conclusion of Theorem~2 hold, under the hypotheses of 
Theorem~2, but replacing the hypothesis on spectra in Theorem~2 by the 
weaker one in Theorem~1?

2) An example of Gramsch \cite{Gr} shows that for the canonical 
homomorphism $\pi : \Cal B(H) \to \Cal C(H)$ (where $H$ is a Hilbert space), 
and for an annulus $U$, it is not always possible to lift an analytic map 
$q : U \to E(\Cal C(H))$ to the whole of~$U$.
Let us suppose in Theorem~1, additionally, that $U$ is simply connected.
Furthermore, both in Theorem~1 and Theorem~2, let us additionally suppose 
that $\pi(\lambda)$ is surjective for each $\lambda \in U$, or for each 
$\lambda \in G$, and, in both cases, let us additionally suppose the 
respective spectral hypotheses for each element 
of\/ $\text{\rm Ker}\, \pi(\lambda)$, for each $\lambda \in U$, 
or for each $\lambda \in G$, respectively.
Can we then choose $V = U$, or $H = G$, respectively?

3) Suppose that all elements in Ker$\, \pi ( \lambda )$, 
for all $\lambda $'s in a set 
which is large enough in some sense, satisfy one of the two topological
spectral hypotheses, from
Theorems 1 and 2. Does it follow, that all elements of all the kernels 
Ker$\, \pi ( \lambda )$, for $\lambda $ in our domain, satisfy the same
spectral hypothesis? 
Observe that in Remark 4 we obtained satisfactory characterizations for
the spectral hypothesis in Theorems 3 and 4. 
The minimal such hypothesis would be for a set which is the union of
sequences of $\lambda $'s in each connected component of our 
domain, converging to points 

\newpage

in the respective connected components of
our domain. For
the hypothesis of quasi-nilpotency, the
sufficiency of this minimal hypothesis, 
for the canonical map to the Calkin
algebra of a Banach space,
was asked in \cite{MZ}. More exactly, let $C( \lambda )$, for $\lambda \in D\,\,
(=$ the open unit disc in ${\Bbb{C}}$), be an
analytic family of {\it{compact operators}} in $\Cal{B}(X)$, for a Banach
space $X$.
If $\sigma \left( C( \lambda ) \right) = \{ 0 \}$ 
for a sequence of $\lambda $'s in $D$, converging to $0$, does the same
equality hold for all $\lambda \in D$?
\cite{MZ} proved that for analytic
families $C( \lambda )$ of {\it{finite rank operators}} the answer is positive.
This might also be related to the problems studied in [AZe]. 

{\bf{Acknowledgements.}} We are grateful to Professors B\'ela Nagy (Szeged), 
Ivan Netuka,
Szil\'ard R\'ev\'esz, and Yuri
Tomilov for useful consultations on potential theory, especially for the
references [AG] and [Ra] cited in Remark 4.

\enddocument